\documentclass{amsart}

\newtheorem{theorem}{Theorem}[section]
\newtheorem{lemma}[theorem]{Lemma}

\theoremstyle{definition}
\newtheorem{definition}[theorem]{Definition}

\newtheorem{thmalp}{Theorem}

\newtheorem{proposition}[theorem]{Proposition}
\theoremstyle{remark}

\numberwithin{equation}{section}

\begin{document}

\title{Brownian motion in a ball in the presence of spherical obstacles}

\author{ Julie O'Donovan}
\address{Department of Mathematics, University College Cork, Cork, Ireland}

\email{j.odonovan@ucc.ie}


\date{June, 2009}


\begin{abstract}
We study the problem of when a Brownian motion in the unit ball
has a positive probability of avoiding a countable collection of
spherical obstacles. We give a necessary and sufficient integral
condition for such a collection to be avoidable.
\end{abstract}

\maketitle

\section{Introduction}
The setting in this paper is the unit ball,  $\mathbb{B}= \{x \in
\mathbb{R}^{d}:|x|<1\}$, in Euclidean space $\mathbb{R}^d$ where
$d \geq 3$. We study the problem of when Brownian motion in the
ball has a positive probability of avoiding a countable collection
of spherical obstacles and thereby reaching the outer boundary of
$\mathbb{B}$.

We denote by $\Lambda$ a sequence of points in
$\mathbb{B}$. To each point $\lambda$ in this sequence we
associate a spherical obstacle, $B(\lambda,r_{\lambda})$, where
\[
B(\lambda, r_{\lambda})=\{x:|\lambda-x|\leq r_{\lambda}\},
\]and denote by
$\partial B(\lambda,r_{\lambda})$ the boundary of this obstacle.
We let $\mathcal{B}$ denote the countable collection of closed
spherical obstacles,
\[
\mathcal{B}=\bigcup_{\lambda \in \Lambda} B(\lambda,r_{\lambda}).
\]
We assume that the spherical obstacles are pairwise disjoint, lie
inside the ball $\mathbb{B}$ and that the origin lies outside
$\mathcal{B}$. We call a collection of spherical obstacles
\emph{avoidable} if there is positive probability that Brownian
motion, starting from the origin, hits the boundary of
$\mathbb{B}$ before hitting any of the spherical obstacles in
$\mathbb{B}$. This is equivalent to positive harmonic measure at 0
of the boundary of the unit ball with respect to the domain
$\Omega=\mathbb{B}\backslash \mathcal{B}$, consisting of the unit
ball less the obstacles, that is $\omega(0,\partial
\mathbb{B};\Omega)>0$.

In the setting of the unit disk, Ortega-Cerd\`{a} and Seip
\cite{CerdaSeip} gave an integral condition for a collection of
disks to be avoidable. In \cite{carrollortega}, Carroll and
Ortega-Cerd\`{a} gave an integral criterion for a configuration of
balls in $\mathbb{R}^d$, $d\geq 3$, to be avoidable.  Thus, it
seems natural to ask if Ortega-Cerd\`{a} and Seip's result for the
disk in the plane can be extended to the ball in space. A solution
to this problem is the main result of this paper.

Next, we put some restrictions on the spacing of the spherical
obstacles. A collection of obstacles, $\mathcal{B}$, is
\emph{regularly spaced} if it is separated, in that there exists
$\epsilon >0$ such that given any $\lambda, \lambda ' \in \Lambda$
with $|\lambda|\geq|\lambda '|$,
 then $|\lambda - \lambda '|>\epsilon (1-|\lambda|)$;
 uniformly dense, in that there exists $R$ with $0<R<1$ such that for $x \in
 \mathbb{B}$, the ball $B(x, R(1-|x|))$ contains at least one $\lambda
 \in \Lambda$; and finally
  the radius $r_\lambda = \phi(|\lambda|)$
 where $\phi:[0,1)\rightarrow [0,1)$ is a decreasing function.

Answering a question of Akeroyd in \cite{Akeroyd}, Ortega-Cerd\`{a}
and Seip \cite{CerdaSeip} proved the following theorem.

\begin{thmalp}\label{2dim}
A collection of regularly spaced disks in the unit disk is avoidable
if and only if
\begin{equation*}\label{oc/seip}
\int_{0}^{1} \frac{\,dt}{(1-t)\log((1-t)/\phi(t))} < \infty.
\end{equation*}
\end{thmalp}
This theorem in \cite{CerdaSeip} is expressed in terms of
pseudo-hyperbolic disks.   We extend Theorem \ref{2dim} to the
setting of the unit ball in $\mathbb{R}^d$, $d\geq 3$.
 \begin{theorem}\label{theorem}
The collection of regularly spaced closed spherical obstacles $
\mathcal{B}$ in $\mathbb{B}$ is avoidable if and only if
\begin{equation}\label{main result}
\int_{0}^{1} \frac{\phi(t)^{d-2}}{(1-t)^{d-1}}\,dt < \infty .
\end{equation}
\end{theorem}

We present two proofs of this result. The first proof exploits a
connection between avoidability and minimal thinness, a potential
theoretic measure of the size of a set near a boundary point of a
region. We learnt of this from both the paper of Lundh \cite{Lundh}
and from Professor S.J. Gardiner. We also make use of a Wiener-type
criterion for minimal thinness
 due to Aikawa \cite{Aikawa}.

The second proof is more direct and transparent. It is an
adaptation of Ortega-Cerd\`{a} and Seip's proof of Theorem
\ref{2dim} in \cite{CerdaSeip}, the key difference being that in
higher dimensions we do not have the luxury of conformal mapping.


\section{Avoidable Obstacles and Minimal Thinness}

Following the notation of Lundh \cite{Lundh}, we let
$SH(\mathbb{B})$ denote the class of non-negative superharmonic
functions on the unit ball and let $P_{\tau}$ denote the Poisson
kernel at $\tau \in \partial \mathbb{B}$. For a positive
superharmonic function $h$ on $\mathbb{B}$ the reduced function of
$h$ with respect to a subset $E$ of $\mathbb{B}$ is
\[R^E_h(w)=\inf \{u(w):u \in SH(\mathbb{B}), u(x)\geq h(x), x \in E
\}\] and the regularized reduced function $\widehat{R}^E_h(w)=\lim
\inf_{x\rightarrow w}R_h^E(x)$.
\begin{definition}
A set $E$ is \emph{minimally thin} at $\tau \in \partial
\mathbb{B}$ if there is an $x_0$ in the unit ball such that
$\widehat{R}^E_{P_{\tau}} (x_0) < P_{\tau} (x_0)$.
\end{definition}
A nice account of reduced functions and minimal thinness may be
found in \cite[Page~38 ff]{Doob} or \cite[Chapter~9]{Gardiner}.

\subsection{Avoidability and minimal thinness} Lundh proves the
following result in \cite{Lundh}. We include a brief proof for the
convenience of the reader.
\begin{proposition}\label{lundh}
Let $A$ be a closed subset of $\mathbb{B}$ such that $\mathbb{B}
\backslash A$ contains the origin and is connected. Let
$\mathcal{M}=\{\tau \in
\partial \mathbb{B}:A \text{ is minimally thin at
}\tau\}$.  Then the following are equivalent:
\begin{itemize}
\item $A$ is avoidable,

\item $|\mathcal{M}| >0$,
\end{itemize}
where $|.|$ denotes surface area on the unit ball.
\end{proposition}

\begin{proof}
Noting that
\[
1= \int_{\partial \mathbb{B}} P_{\tau} (x) \frac{\,d\tau}{|\partial
\mathbb{B}|},
\]
and taking $h\equiv1$ in \cite[Corollary 9.1.4]{Gardiner}, we see
that
\[ \widehat{R}_{1}^{A} (x)=
 \int_{\partial \mathbb{B}} \widehat{R}_{P_{\tau}}^{A} (x)
 \,\frac{d\tau}{|\partial \mathbb{B}|}.
\] Also, it follows from \cite[Page 653, 14.3sm]{Doob} that the regularized reduced function
of $1$ with respect to $A$ evaluated at $x$ is the harmonic
measure at $x$ of $\partial A$ in the domain $\mathbb{B}
\backslash A$. Thus,
\begin{equation*}
\omega(0,\partial A,\mathbb{B} \backslash A)=\widehat{R}_{1}^{A}
(0)=
 \frac{1}{|\partial \mathbb{B}|} \int_{\partial \mathbb{B}} \widehat{R}_{P_{\tau}}^{A} (0)
 \,d\tau.
\end{equation*}
Since $\widehat{R}_{P_{\tau}}^{A} (0)\leq P_{\tau}(0)=1$, it follows
that $\omega(0,\partial A,\mathbb{B} \backslash A)<1$ if and only if
the set $\mathcal{M}_0=\{\tau \in
\partial \mathbb{B},\widehat{R}_{P_{\tau}}^{A}(0)<1\}$
has positive measure. In the connected domain $\mathbb{B}\backslash
A$, the set $\mathcal{M}_0$ is the same as the set $\mathcal{M}$.
Thus, $A$ being avoidable, that is $\omega(0,\partial
\mathbb{B};\mathbb{B}\backslash A)>0$, is equivalent to
$\mathcal{M}$ having positive measure.
\end{proof}

\subsection{Minimal thinness and a Wiener-type criterion}
 It is a standard result, see for example Aikawa \cite{Aikawa}
or Lundh \cite{Lundh}, that a set is minimally thin at a point if
and only if it satisfies a Wiener-type criterion. Let $\{Q_k\}$ be
a Whitney decomposition of the unit ball $\mathbb{B}$ in
$\mathbb{R}^d$ ($d\geq 3$) and let $q_k$ be the Euclidean distance
from the centre, $c_k$, of the Whitney cube $Q_k$ to the boundary
of $\mathbb{B}$. Let $A$ be a subset of $\mathbb{B}$. Let $\tau$
be a boundary point of $\mathbb{B}$ and $\rho_k (\tau)$ be the
distance from $c_k$ to the boundary point $\tau$. Let cap denote
Newtonian capacity. Then $A$ is minimally thin at the point $\tau$
if and only if
\begin{equation}\label{WC}
\sum_{k} \frac{q_k ^2}{\rho_k(\tau)^d} \text{cap} (A \cap
Q_k)<\infty.
\end{equation}
 In the next section, we consider this Wiener-type criterion in the
particular setting of the unit ball less a collection of regularly
spaced spherical obstacles.
\subsection{ Wiener-type criterion and integral condition  }
For a constant $K>1$, we let $S_{j}=\{x:|x|=1-K^{-j}\}$ be the
sphere of radius $\rho_j=1-K^{-j}$ and $B_{j}$ be the interior of
this sphere. We denote by $A_j$ the annulus bounded by $S_{j}$ and
$S_{j-1}$, and write $\phi_j$ for $\phi(\rho_j)$.
\begin{proposition} \label{WC and ic}
Let $\mathcal{B}$ be a regularly spaced collection of spherical
obstacles in $\mathbb{B}$.  \\
(i) If the set $\mathcal {B}$ satisfies the Wiener-type criterion
(\ref{WC}) at some point in $\partial \mathbb{B}$ then the integral
condition (\ref{main result}) holds,
\\  (ii) The integral condition (\ref{main result}) implies that $\mathcal {B}$ satisfies the
Wiener-type criterion (\ref{WC}) at all points $\tau \in \partial
\mathbb{B}$.
\end{proposition}

\begin{proof}
We first assume that the integral condition holds and we'll show
that (\ref{WC}) follows. We note that the integral condition
(\ref{main result}) is equivalent to
\begin{align}\label{int-sum}
\sum _{j=1}^{\infty} (\phi_j K^{j})^{d-2} < \infty,
\end{align} where $K>1$. By the separation condition on the
sequence $\Lambda$, there is an $N$ such that any cube $Q_k$ can
contain no more than $N$ points in $\Lambda$.
 Splitting the sum in (\ref{WC}) into a sum over annuli we obtain
\begin{align}
\sum_{k} \frac{q_k ^2}{\rho_k(\tau)^d} \text{cap} (\mathcal{B}
\cap Q_k)&= \sum_{j=1}^{\infty} \sum_{k:c_k\in A_j} \frac{q_k
^2}{\rho_k(\tau)^d} \text{cap} (\mathcal{B} \cap Q_k)
\label{wcseries}
\\& \leq \sum_{j=1}^{\infty}N(K^{-j})^2 \phi_j^{d-2}  \sum_{k: c_k\in
A_j} \frac{1}{\rho_k(\tau)^d}, \label{lattersum}
\end{align}
 since the capacity of a ball with radius
$\phi_j$ is equal to $\phi_j^{d-2}$. We now concentrate on the
latter sum in (\ref{lattersum}). We split up the $j^{th}$ annulus
$A_j$ into rings centered at the projection of $\tau$ onto the
sphere $S_j$, and with radius equal to $n K^{-j}$ where we recall
that $K^{-j}$ is the distance from $\tau$ to $S_j$. There are at
most
\[\frac{c_d (n K^{-j})^{d-2}}{(K^{-j})^{d-2}}=c_dn^{d-2}\]
Whitney cubes in each ring where $c_d$ is a constant depending on
the dimension, $d$. For the $n^{th}$ ring,
 \[\rho_k (\tau) \geq
nK^{-j}\] and $N_j$ rings intersect the annulus $A_j$. Thus,
\begin{align*}
\sum_{k: c_k\in A_j} \frac{1}{\rho_k (\tau)^d} \leq&
\sum_{n=1}^{N_j}
\frac{c_d n^{d-2}}{(nK^{-j})^d} \\
\leq & (K^j)^d c_d\sum_{n=1}^{N_j} \frac{1}{n^2}.
\end{align*}
Thus, we see that the Wiener-type series (\ref{wcseries}) is
convergent.

We now assume that the set $\mathcal{B}$ satisfies (\ref{WC}) at
some arbitrary point $\tau \in \partial \mathbb{B}$ and show that
this implies the integral condition (\ref{main result}). We choose
$K$ sufficiently large so that for all $j$ bigger than a fixed
constant there is at least one centre of a ball in each Whitney
cube, $Q_k$, in the resulting Whitney decomposition of
$\mathbb{B}$. Starting with the Wiener-type series we
 split it into a sum over the annuli $A_j$ and then proceed to
 ignore all Whitney cubes in $A_j$ except one near to the
point $\tau$, for which $\rho_k(\tau) \leq K^{-j}$, as follows.
\begin{align*}
\sum_{k} \frac{q_k ^2}{\rho_k(\tau)^d} \text{cap} (\mathcal{B}
\cap Q_k)& = \sum_{j=1}^{\infty} \sum_{k: c_k\in A_j} \frac{q_k
^2}{\rho_k(\tau)^d} \text{cap} (\mathcal{B} \cap Q_k) \\
&\geq \sum_{j=0}^{\infty}K^{-2j} \phi_j^{d-2}
\frac{1}{\rho_k(\tau)^d}\\
&\geq\sum_{j=0}^{\infty} (\phi_j K^{j})^{d-2}
\end{align*}
Thus, since the Wiener-type series is convergent, (\ref{int-sum})
follows and so the integral condition (\ref{main result}) holds.
\end{proof}
Combining Proposition \ref{lundh}, the Wiener-type criterion
(\ref{WC}) and Proposition \ref{WC and ic} we have a proof of
Theorem \ref{theorem}. We note that the method used in this
section could also be used to give an alternative proof of
Ortega-Cerd\`{a} and Seip's Theorem~\ref{2dim}.


\section{Direct proof of Theorem \ref{main result}}
We now give an alternative proof of Theorem \ref{theorem} by
adapting the method of Ortega-Cerd\`{a} and Seip  in
\cite{CerdaSeip}. In dimensions higher than 2 we do not have
conformal mapping, but we do have the Kelvin transform. We let
\[x^*= \frac{\rho_{j+1}^2}{|x|^2}x\] be the inversion of the point
$x$ in the sphere of radius $\rho_{j+1}$. We note that $|x||x^*|$
equals $\rho_{j+1}^2$, and let $\phi(|\lambda|)=\phi_{\lambda}$. We
begin with a lemma, prove the sufficiency of the integral condition
in the next subsection and the necessity in the following one.

\begin{lemma}\label{lamda-lamba*}
Let $K>\max\{4,\frac{1+R}{1-R}\}$ and $x$ be an arbitrary point
belonging to $S_{j-1}$. There is a centre of an obstacle,
$\lambda_x \in \Lambda$, such that $\lambda_x$ lies in the annulus
$A_j$ bounded by $S_{j-1}$ and $S_j$, and
\[|x-\lambda_x|\leq \frac{K-1}{K}|x^*-\lambda_x|.\]
\end{lemma}

\begin{proof}
For $x \in S_{j-1}$, let $x'$ be the point on the extension of the
radius of $S_j$ containing $x$, and located halfway between
$S_{j-1}$ and $S_{j}$. Then $x'$ is a distance
$K^{-(j-1)}-\frac{K-1}{2K^j}$ from the boundary of the ball
$\mathbb{B}$. Since $\Lambda$ is uniformly dense, the ball
$B(x',R(1-|x'|))$ contains some $\lambda_x \in \Lambda$. Also, due
to the choice of $K$,
 the ball $B(x',R(1-|x'|))$ is contained in the annulus $A_j$.
Let $x''$ be on the same ray as $x$ and $x^*$ and also on $S_j$.
We first note that $|x-\lambda_x|\leq|x-x''|$ and
$|x^*-\lambda_x|>|x^*-x''|$. Also, we note that $|x|=\rho_{j-1}$,
$|x''|=\rho_j$ and $|x^*|=\rho_{j+1
 }^2/\rho_{j-1}$. Thus,
\begin{equation}\label{near ball}
|x-\lambda_x|\leq|x-x''|= (K-1)K^{-j}. \end{equation}
 Also,
\begin{align*}
|x^*-\lambda_x|\geq|x^*-x''|&= \frac{(1-K^{-(j+1)})^2}{1-K^{-(j-1)}}
- (1-K^{-j}) \geq  K^{-j+1},
\end{align*}
for
 $j\geq2$.
 Thus,
 \[|x-\lambda_x|\leq
\frac{K-1}{K}|x^*-\lambda_x|,\] as required.
\end{proof}

\subsection{Integral Condition (\ref{main result}) implies avoidability}

We first assume (\ref{main result}) and show that the spherical
obstacles are avoidable that is, we show that $\omega(0,
\partial \mathbb{B}; \Omega)
>0 $.  We split the collection of spherical obstacles into
those with centres inside and those with centres outside a ball of
radius $r<1$. We let $\Lambda_{r} = \{\lambda\in \Lambda:
|\lambda|>r \}$ and let \[\mathcal{B}_r=\bigcup_{\lambda \in
\Lambda_{r}}B(\lambda, r_{\lambda})=\bigcup_{\lambda \in
\Lambda_{r}}B_{\lambda} \] denote the infinitely many spherical
obstacles with centres outside $B(0,r)$. Also, we let
$\Omega_{r}=\mathbb{B} \backslash \mathcal{B}_r$ be the champagne
subregion where all obstacles have centres outside a ball of
radius $r$. We may safely ignore the finitely many spherical
obstacles with centres inside the ball of radius $r$. Thus, it is
sufficient to show that $\omega(0,
\partial \mathbb{B}; \Omega_{r})> 0$ for some $r$ with $0<r<1$, which
is equivalent to showing that $\omega(0,
\partial  \mathcal{B}_r; \Omega_{r})<1$. We choose $r$ such that \[\int_{r}^{1}
\frac{\phantom{{x}}\phi(t)^{d-2}}{(1-t)^{d-1}}dt<
\frac{\epsilon^{d}(K-1)^{d-2}}{2^{d+1}d(d-2)K^{2d-1}}\] and let
$n_r$ be the biggest integer smaller than
 $1+\log(\frac{1}{1-r})/\log K$. This ensures that $r>\lfloor 1-K^{-(n_r-1)}\rfloor$.
We proceed as follows,
\begin{align*}
\omega(0, \partial \mathcal{B}_r ; \Omega_r) = \sum_{\lambda \in
\Lambda_r} \omega (0, \partial B_{\lambda};\Omega_r) \leq
&\sum_{\lambda \in \Lambda_r} \omega (0,\partial
B_{\lambda};\mathbb{B}
\backslash B_{\lambda}) \\
\leq& \sum_{j=n_r}^{\infty}( \sum_{\lambda \in A_j}\omega (0,
\partial B_{\lambda}; \mathbb{B}\backslash B_{\lambda})).
\end{align*}
We now obtain an upper bound for the number of centres in $A_j$
and an upper bound for the contribution of an obstacle with centre
in $A_j$ to the above sum. Due to the separation condition,
centres of balls in $A_j$ are at least $\epsilon K^{-j}$ apart.
Thus, the number of centres in $A_j$, which is less than the
volume of $A_j$ divided by the volume of a ball with radius
$\epsilon K^{-j}/2$, is less than
\[\frac{2^d d K^{2}}{\epsilon^d}K^{(d-1)j}.\]
Next, we want an upper bound for $\omega(0,\partial
B_{\lambda};\mathbb{B}\backslash B_{\lambda})$. We construct a
suitable function $h$ that is harmonic on $\mathbb{B} \backslash
B_{\lambda}$, continuous on its closure and also satisfies $
h(x)\geq 1$, $x \in
\partial B_{\lambda}$   and $ h(x) \geq 0$, $x \in \partial
\mathbb{B}$. Then, using the Maximum Principle, we obtain the
required upper bound. Consider the function
 \[h(x) = 2\left[ u_\lambda(x) -
u^* _\lambda (x)\right],\] where
\[u_\lambda(x) =
\left[\frac{\phi_{\lambda}}{|x-\lambda|}\right]^{d-2}\text{, }
u^*_\lambda(x) =
\left[\frac{\phi_{\lambda}}{|x||x^*-\lambda|}\right]^{d-2} \text{
and } x^*= \frac{1}{|x|^2}x.\] We note that $u_\lambda$ and
$u^*_\lambda$ are harmonic. Also, $1/2$ is a lower bound for $
u_\lambda(x) - u^* _\lambda (x)$ for $x \in
\partial B_\lambda$ which we show as follows. For $x \in
\partial B_{\lambda}$, we have that $|x|\geq 1-K^{-1}$
 and $|x^*-\lambda|\geq K^{-j}$, hence
\begin{align*}
u_\lambda(x) - u^* _\lambda (x) &=
1-\left[\frac{\phi_{\lambda}}{|x||x^*-\lambda|}\right]^{d-2}\geq
1-\left[\frac{K\phi_{j-1}}{(K-1)K^{-j}}\right]^{d-2}.
\end{align*}
It follows from (\ref{int-sum}) that
\[\lim_{j\rightarrow \infty} \frac{\phi_{j-1}}{K^{-j}}=0.\]
Thus, there exists $N$ such that for $j>N$
\begin{align*}
u_\lambda(x) - u^* _\lambda (x) >\frac{1}{2}.
\end{align*}
Thus, $h(x)$ satisfies the required criteria and is an upper bound
for the harmonic measure $\omega(x,
\partial B_{\lambda};\mathbb{B}\backslash B_{\lambda})$.

 Next, we want an
upper estimate for $h(0)$. We first note that as $x \rightarrow 0$,
$x^* \rightarrow \infty$ and also that $|x||x^*|=1$. Thus, as $x
\rightarrow 0$, $u^* _\lambda (x) \rightarrow \phi_{\lambda}^{d-2}$.
Next,
\begin{align*}
\frac{1}{2}h(0)= & \lim_{x\rightarrow 0}[u_\lambda(0) - u^*
_\lambda
 (0)]=\left(\frac{\phi_{\lambda}}{|\lambda|}\right)^{d-2}
-\phi_{\lambda}^{d-2} =
\phi_{\lambda}^{d-2}\left[\frac{1-|\lambda|^{d-2}}{|\lambda|^{d-2}}\right]\\\leq
&
\left(\frac{\phi_{j-1}}{|\lambda|}\right)^{d-2}(d-2)\left[K^{-(j-1)}+O(K^{-2j})\right].
\end{align*}
Thus, for sufficiently large $j$,
\[h(0)\leq
4K(d-2)\left(\frac{\phi_{j-1}}{1-K^{-(j-1)}}\right)^{d-2}K^{-j}.\]
Therefore,
\begin{align*}
\omega(0, \partial \mathcal{B}_r ; \Omega_r) & \leq
\sum_{j=n_r}^{\infty}\frac{2^d d K^{2}}{\epsilon^d}K^{(d-1)j}
 4K(d-2)\left(\frac{\phi_{j-1}}{1-K^{-(j-1)}}\right)^{d-2}K^{-j}
 \\&\leq \frac{2^{d+2}d(d-2)K^{2d-1}}{\epsilon^{d}(K-1)^{d-2}}
 \sum_{j=n_r}^{\infty}\left(\phi_{j-1} K^{j-1}\right)^{d-2} <1
\end{align*}
provided $n_r$ is suitably selected as described at the start of the
proof. Thus, $\omega(0,\partial \mathcal{B}_r;\Omega_{r})<1$ and
hence we see that $\omega(0, \partial \mathbb{B};\Omega)>0$ as
required.

\subsection{Avoidability implies the integral condition (\ref{main result})}

 Now we assume that $\omega(0,\partial\mathbb{B};\Omega)>0$
 and we'll show (\ref{main result}) holds. We begin by ignoring all obstacles with centres in an annulus
  $A_j$ where $j$ is
odd. We let
 \[\Omega'=\mathbb{B}\backslash \bigcup_{\lambda \in A_j , \text{ j even}}B(\lambda,
 r_\lambda)\]and note that since
$\omega(0,
\partial \mathbb{B}; \Omega)>0$, then $\omega(0,
\partial \mathbb{B}; \Omega')>0$.  We choose $K>\max\{4,\frac{1+R}{1-R}\}$,
 where $R$ is the constant mentioned in the definition of
regularly spaced.
 We let
 $P_j$ denote the probability that Brownian motion starting at the
 origin hits $S_{j+1}$ before hitting any of the obstacles with centres in
 $B_j$ but not in any $A_i$ where $i$ is odd. We let $Q_j$ denote the supremum of the probabilities that Brownian
 motion starting on $S_{j-1}$ hits $S_{j+1}$ before hitting any of the
obstacles with centres in $A_j$.  We note that $P_j \leq Q_j
P_{j-2}$ and
 that therefore for $n$ even
\[P_n \leq  P_0 \prod_{j=1, \text{ j even}}^{n} Q_j.\] Since $\omega(0,
\partial \mathbb{B}; \Omega')=\delta>0$, it follows that $P_n\geq\delta$ for all $n$ and,
since $Q_j < 1$,
 \begin{equation}\label{even sum finite}
 \sum_{j=1, \text{ j even}}^{\infty} (1-Q_j) < \infty.
 \end{equation}
  We note
that $1-Q_j$ is the infimum over $x \in S_{j-1}$ of the probability
that Brownian motion starting at $x$ hits a ball with centre in
$A_j$ before hitting $S_{j+1}$. Thus, if we consider only a single
ball near $x$, say $B_{\lambda_x}$ where $\lambda_x$ is the centre
of the ball near $x$ as described in Lemma \ref{lamda-lamba*},  then
\[1-Q_j \geq \inf_{x \in S_{j-1}} \omega(x, \partial B_{\lambda_x};
 B_{j+1} \backslash B_{\lambda_x}).\] Thus, we need a lower
bound for $\omega(x,
\partial B_{\lambda_x}; B_{j+1} \backslash B_{\lambda_x})$.
 We want a suitable function $h_j$ that is harmonic on $B_{j+1}
\backslash B_{\lambda_x}$, continuous on its closure and also
satisfies $h_j(y)\leq 1$, $y \in \partial B_{\lambda_x}$ and $h_j(y)
\leq 0$, $y \in S_{j+1}$. Then we can again avail of the Maximum
Principle to obtain the required lower bound. Consider the function
 \[h_j(y) =  u_\lambda(y) -
u^* _\lambda (y),\] where
\[u_\lambda(y) =
\left[\frac{\phi_{\lambda}}{|y-\lambda_x|}\right]^{d-2}\text{, }
u^*_\lambda(y) =
\left(\frac{\rho_{j+1}}{|y|}\right)^{d-2}\left[\frac{\phi_{\lambda}}{|y^*-\lambda_x|}\right]^{d-2}
\text{ and } y^*= \frac{\rho_{j+1}^2}{|y|^2}y.\]
 Then $h_j(y)$ satisfies the required criteria as both $u_{\lambda}$
 and $u^*_\lambda$ are harmonic, $h_j\leq u_{\lambda}=1$ on
 $\partial B_{\lambda_x}$, and $u_{\lambda}=u^*_\lambda$ on
 $S_{j+1}$. Next, we want a lower estimate for $h_j$ at the point $x\in
S_{j-1}$. With the help of Lemma \ref{lamda-lamba*},
\begin{align*}
 u_{\lambda}(x)-u_{\lambda}^*(x) & =
 \left[\frac{\phi_{\lambda}}{|x-\lambda_x|}\right]^{d-2}-
\left(\frac{\rho_{j+1}}{\rho_{j-1}}\right)^{d-2}\left[\frac{\phi_{\lambda}}{|x^*-\lambda_x|}\right]^{d-2} \\
& \geq \left(\frac{\phi_j}{|x-\lambda_x|}\right)^{d-2}\left[1-
\left(\frac{\rho_{j+1}}{D \rho_{j-1}}\right)^{d-2}\right],
\end{align*}
where $D=K/(K-1)>1$. Then for sufficiently large $j$, namely $j$
where
\[\frac{\rho_{j+1}}{\rho_{j-1}}<\frac{1+D}{2},\]
we find that
\[ u_{\lambda}(x)-u_{\lambda}^*(x)  \geq
c\left(\frac{\phi_j}{|x-\lambda_x|}\right)^{d-2},\]where $c$ is
some
 positive constant.

By (\ref{near ball}), we find that for $x \in S_{j-1}$,
\begin{align*}
\omega(x,\partial B_{\lambda_x};B_{j+1}\backslash B_{\lambda_x})
\geq h_j(x) = u_{\lambda}(x)-u_{\lambda}^*(x)\geq c(K-1)^{2-d}
(\phi_j K^{j})^{d-2}.
\end{align*}
It now follows from (\ref{even sum finite}) that
\[\sum_{j=1, \text{ j even}}^{\infty} (\phi_j K^{j})^{d-2} < \infty.\]
Similarly it may be shown that  \[\sum_{j=1, \text{ j
odd}}^{\infty} (\phi_j K^{j})^{d-2} < \infty,\] and so
 \[\sum_{j=1}^{\infty}
(\phi_j K^{j})^{d-2} < \infty.\] Hence, (\ref{main result}) holds
and the proof is complete.

\section*{Acknowledgements}
The author would like to thank Tom Carroll for his help and
guidance throughout this work.


\end{document}